\documentclass[10pt]{amsart}
\usepackage{amsmath,amssymb,amsxtra,amsthm,amscd}
\usepackage{lmodern}
\usepackage[all,cmtip]{xy}
\usepackage{rotating}
\usepackage{microtype}
\newif\ifShowLabels
\ShowLabelstrue
\newcommand{\TeXref}[1]{
\marginpar{\scriptsize \texttt{#1}}}

\DeclareMathOperator{\diam}{diam}

\DeclareMathOperator{\K}{\mathit{K}}

\DeclareMathOperator{\V}{\mathcal{V}}

\DeclareMathOperator*{\one}{1}
\newcommand{\onehatplace}[1]
{ \one^{\substack{#1 \\ \frown}} }

\DeclareMathOperator*{\bones}{\times}
\newcommand{\undertimes}[1]
{ \bones_{#1} }

\DeclareMathOperator*{\bowl}{\cup}
\newcommand{\undercup}[1]
{ \bowl_{#1} }

\DeclareMathOperator*{\arch}{\cap}
\newcommand{\undercap}[1]
{ \arch_{#1} }

\newcommand{\pull}
{\!\!\! -\!\!\! -\!\!\! -\!\!\!}

\DeclareMathOperator*{\holimprep}{holim}                       
\newcommand{\holim}[1]%
{\displaystyle\holimprep_{\substack{\leftarrow \pull - \\ #1}} \, }

\DeclareMathOperator*{\hocolimprep}{hocolim}                   
\newcommand{\hocolim}[1]%
{\displaystyle\hocolimprep_{\substack{- \pull \rightarrow \\ #1}} \, }

\DeclareMathOperator*{\plainlim}{lim}                           
\newcommand{\contralim}[1]%
{\displaystyle\plainlim_{\substack{\leftarrow \pull - \\ #1}} \, }

\DeclareMathOperator*{\plaincolim}{colim}                       
\newcommand{\colim}[1]%
{\displaystyle\plaincolim_{\substack{- \pull \rightarrow \\ #1}} \, }

\DeclareMathOperator*{\laxlimplain}{laxlim}                     
\newcommand{\laxlim}[1]%
{\displaystyle\laxlimplain_{\substack{\leftarrow \pull - \\ #1}} \, }

\theoremstyle{plain}
\newtheorem{Thm}{Theorem}[section]

\newtheorem{Cor}[Thm]{Corollary}

\newtheorem{Lem}[Thm]{Lemma}
\newtheorem{Prop}[Thm]{Proposition}

\theoremstyle{definition}
\newtheorem{Def}[Thm]{Definition}

\newtheorem{Ex}[Thm]{Example}

\newtheorem{Rem}[Thm]{Remark}

\theoremstyle{remark}
\newtheorem{Not}[Thm]{Notation}

\newtheoremstyle{freestylethm}{6pt}{6pt}{\itshape}{}%
                {\bfseries}{}{.5em}{\thmnote{#3}}
\theoremstyle{freestylethm}

\newcommand{\SecRef}[2]{\section{#1}\label{S:#2}%
\ifShowLabels \TeXref{{S:#2}} \fi}
\newcommand{\SSecRef}[2]{\subsection{#1}\label{SS:#2}%
\ifShowLabels \TeXref{{SS:#2}} \fi}

\newcommand{\refS}[1]{\textup{\ref{S:#1}}}
\newcommand{\refSS}[1]{\textup{\ref{SS:#1}}}

\newcommand{\refT}[1]{\textup{\ref{T:#1}}}

{ \begin{Thm} \label{T:#1}
\ifShowLabels \TeXref{T:#1} \fi }%
{ \end{Thm} }
\newenvironment{DefRef}[1]%
{ \begin{Def} \label{D:#1}
\ifShowLabels \TeXref{D:#1} \fi }%
{ \end{Def} }
{ \begin{Lem} \label{L:#1}
\ifShowLabels \TeXref{L:#1} \fi }%
{ \end{Lem} }
{ \begin{Cor} \label{C:#1}
\ifShowLabels \TeXref{C:#1} \fi }%
{ \end{Cor} }
\newenvironment{RemRef}[1]%
{ \begin{Rem} \label{R:#1}
\ifShowLabels \TeXref{R:#1} \fi }%
{ \end{Rem} }
{ \begin{Prop} \label{P:#1}
\ifShowLabels \TeXref{P:#1} \fi }%
{ \end{Prop} }
{ \begin{Ex} \label{E:#1}
\ifShowLabels \TeXref{E:#1} \fi  }%
{ \end{Ex} }
\newenvironment{NotRef}[1]%
{ \begin{Not} \label{N:#1}
\ifShowLabels \TeXref{N:#1} \fi }%
{ \end{Not} }

\newenvironment{ThmRefName}[2]%
{ \begin{Thm} [#2]\label{T:#1}
\ifShowLabels \TeXref{T:#1} \fi }%
{ \end{Thm} }
\newenvironment{DefRefName}[2]%
{ \begin{Def} [#2]\label{D:#1}
\ifShowLabels \TeXref{D:#1} \fi }%
{ \end{Def} }
{ \begin{Lem} [#2]\label{L:#1}
\ifShowLabels \TeXref{L:#1} \fi }%
{ \end{Lem} }
{ \begin{Cor} [#2]\label{C:#1}
\ifShowLabels \TeXref{C:#1} \fi }%
{ \end{Cor} }
{ \begin{Rem} [#2]\label{R:#1}
\ifShowLabels \TeXref{R:#1} \fi }%
{ \end{Rem} }
{ \begin{Prop} [#2]\label{P:#1}
\ifShowLabels \TeXref{P:#1} \fi }%
{ \end{Prop} }
{ \begin{Ex} [#2]\label{E:#1}
\ifShowLabels \TeXref{E:#1} \fi }%
{ \end{Ex} }

 \ShowLabelsfalse

\begin{document}

\title[Topological methods in robot motion path planning]{Some geometric and topological data-driven methods in robot motion path planning}
\author[Boris Goldfarb]{Boris Goldfarb}
\address{Department of Mathematics and Statistics\\ State University of New York\\ Albany\\ NY 12222-0001}
\email{bgoldfarb@albany.edu}
\date{\today}

\begin{abstract}
Motion path planning is an intrinsically geometric problem which is central for design of robot systems.  Since the early years of AI, robotics together with computer vision have been the areas of computer science that drove its development.  Many questions that arise, such as existence, optimality, and diversity of motion paths in the configuration space that describes feasible robot configurations, are of topological nature.  The recent advances in topological data analysis and related metric geometry, topology and combinatorics have provided new tools to address these engineering tasks.  We will survey some questions, issues, recent work and promising directions in data-driven geometric and topological methods with some emphasis on the use of discrete Morse theory.
\end{abstract}

\maketitle

\SecRef{Introduction}{intro}

We start by establishing the problem of motion path planning.  There are two types of paths one considers in robotics.  The notion of path most closely associated with real tasks is the \textit{underactualized} notion where the path is viewed as a continuous flow of instructions that control the internal degrees of robot joints.  This can be viewed as the final product of our analysis prescribing a motion algorithm.  In general, this algorithm doesn't directly translate to the coordinates on a path in the parameter space of all possible unobstructed configurations of the robot. One good example is a snake robot.  The underactualized path describes the coordinated motions of all joints but the actual modeled movement depends on other physical parameters of the environment such as the friction characteristics of the skin and the floor.  In lots of situations however the more direct and generally simpler problem of finding a path in the parameter space of the robot is essentially the same as defining a motion algorithm.  This is true for moving robot arms bolted to a factory floor or even, with some interpretation, self-driving car control instructions.  We will survey methods that study the simpler, \textit{actualized} path search problem, which is evidently a purely geometric problem.  

Even in this simplified context there is a variety of important complications that we will not address.  While it will be elementary to prescribe structural restrictions that may be imposed, some of the parameters may need to be velocities of the robot components, and so the imposed obstructions may be constraints on the magnitude of these velocities.  The parameter space in this case is often modeled on a subspace of the tangent bundle on the genuine space of structural positions.  Another complication can be the need to model several robots, which can be of different design, moving in the same physical environment at the same time.

The basic case we consider is one single robot whose configuration is described by giving $n$ numerical parameters $x = (x_1$, $\ldots$, $x_n)$ in the $n$-dimensional Euclidean space $E$.  Obstacles in this setting is the set $\mathcal{O}$ of all configurations that are not allowed because they are either obstructed by obstacles in the robot environment, not achievable because of the structure of the robot itself due to self-obstructions, or just postulated to be not allowed.  

The \textit{free space} or \textit{free configuration space} $X = E \setminus \mathcal{O}$ is usually called \textit{C-space}.  This is where path planning is done.  We will need to introduce a metric in C-space.  Let's assume that the metric in $X$ is the subspace metric from the Euclidean space.  The methods we will concentrate on are quite agnostic to the metric used, so the choice can be made tailored to the needs of a specific task.  

Sometimes what we are after is described as a ``piano movers problem''.  This supplies intuition which makes clear that small continuous deformations of the obstacles can lead to discontinuous---often extremely jarring in simulations---changes in C-space, even in the homotopy type of the C-space.  

\subsection*{Acknowledgements and outline}
This author's personal interest in this subject developed from collaborations \cite{UGE:21,UGE:22,UGWE:20,UWE:19}, where discrete Morse theory was introduced in motion path planning.  We have demonstrated this use to be competitive in practice against the state-of-the-art methods.  I am also convinced this application of discrete Morse theory can be made very systematic.  I am grateful to my collaborators for this experience.  I am grateful to Matt Zaremsky for conversations about his paper \cite{mZ:22}, related to open problem (9) in Section \refS{problems}.  I am grateful to the referee for suggesting improvements to the paper.

Before surveying the applications of Morse theory in Section 4, 
I will indicate the variety of methods of topological origin used in robotics in Section 2.  
Section 3 contains the necessary background from discrete Morse theory and some general comments on its relevance.
Section 5 collects problems and questions that are at the heart of motion path planning in robotics and are intrinsically geometric.

\SecRef{Panorama of topological methods in robot motion planning}{Pano}

\SSecRef{Topological complexity}{TC}
One major body of work is due to Michael Farber, Shmuel Weinberger, Jesus Gonzales, their collaborators and other topologists who study topological complexity of spaces.  The main focus is to learn to distribute path planning in complicated spaces by covering the space with contractible patches where planning has essentially a unique solution that depends continuously on the initial and terminal configurations.  The main thrust of this research is to learn the minimal number of patches required for a given space.  There are several papers in this collection that review this work and various related developments, so we are compelled to move on to other topics complementary to this large area.

\SSecRef{Artificial potential vector fields}{PVF} 
Mainstream robotics researchers have been exploiting topology and geometry in the form of smooth potential flows on a smooth C-space for a very long time.  Among the most active and prominent proponents of geometric methods are Howie Choset, Daniel Koditschek, Steven LaValle.  Koditschek \cite{K:89b} traces the idea of introducing artificial potential flows to induce directed motion while avoiding obstacles to the early work of Khatib \cite{oK:78,oK:86}.  All modern textbooks about motion planning have a central chapter on these ideas, cf. chapter 8 in LaValle \cite{sL}, chapter 4 in Choset et al \cite{hC:05}, chapter 10 in Lynch/Park \cite{kLfP}.  

The smooth methods such as Koditschek et al. \cite{K:89a,K:89b,KR:90,KR:91,KR:92} are widely applied, with citation counts for flagship papers in several thousands.  There are also extensions of the use of the navigation Morse functions to create so-called Morse decompositions of the C-space as in Acar et al. \cite{eA:02}.  The basic idea is to create a function on C-space that would mimic the potential scalar field in space.  The goal is to derive the corresponding differentiable vector field with a source at the initial point (say, where the potential has a global maximum) and a sink at the terminal point of the motion (where the potential has the unique global minimum).  The trajectories along the flow lines of the field should then give a solution to the motion path problem.  The artificial nature of the vector field allows one to arrange various desired features.  For example, declaring potential values that are or near maximal on the boundaries of the obstacles guarantees avoidance of obstacles, and even better, allows to control the proximity of the flow lines to the obstacles.  This latter arrangement becomes important in questions of safety of the designed motion so a certain degree of separation of trajectories from the obstacles is maintained.  The choice of values on the obstacles allows to modulate between optimality and safety of the motion paths.  The central theorem can be stated as follows.

\begin{ThmRefName}{KThm}{Koditschek \cite{K:89a}, Koditschek/Rimon \cite{KR:90,KR:91}}
    Let $M \subset E^n$ be a compact connected $n$-dimensional analytic submanifold with boundary. Suppose we are given a navigation function, that is, an analytic Morse function $\phi \colon M \to [0, 1]$ with a unique minimum in the interior point $p$ of $M$, and with the maximum value attained on all boundary components.  Then the control laws resulting from $\phi$ define a closed loop robotic system whose trajectories approach the destination $p$ without intersecting obstacles, starting in an open dense set of initial states.  Moreover, if $h \colon N \to M$ is an analytic diffeomorphism then the pullback $\phi \circ h \colon N \to [0, 1]$ is a navigation function on $N$.
\end{ThmRefName}

To the best of my knowledge, this method remains the only provably correct general method for motion path planning.  The challenge in practice is to create a navigation function with the required properties.  One approach commonly used is to start with a partially defined function and adjust as one observes a trajectory encountering a local minimum which is not the global minimum $p$.  There exist paradigms that create a parameterized family of functions that can be used to avoid or get out of a local minimum.

\SSecRef{Topological data analysis}{TDA} 
Finally, there is a variety of data-driven methods that are usually based on a sample drawn from C-space.  

For simplicity, let's still assume C-space is a subspace of the Euclidean space $E^n$ with the restricted Euclidean metric.  

\begin{NotRef}{cubulation}
If the C-space is a codimension $0$ submanifold $M$, it's common in robotics to generate a uniform sample in $M$ by considering a fine $n$-dimensional grid with a fixed length of the edges $l > 0$ in $E^n$ and intersecting its vertex set with $M$.  We will denote this discrete sample $D_{\mathrm{cube},l}$.  Instead of generating a discrete sample from the vertices, one can also consider the set of complete $n$-cubes from the corresponding cubulation of $E^n$ that are contained entirely within $M$.  Then $M$ is approximated by an $n$-dimensional cubical complex denoted $C_{M,l}$.  
\end{NotRef}

Other options are from a variety of random sampling methods used to build a so-called map of $M$. The map itself is usually presented in the form of a graph where the vertices are a sample from $M$.  The intrinsic methods would usually require a known specific geometry of $M$ and so are not of interest to us.  We want no assumption on $M$. 

The two main types of sampling-based techniques are the Probabilistic RoadMap (PRM) and the Rapidly-exploring Random Tree (RRT) algorithms. Section 10.5 in Lynch--Park \cite{kLfP} is a good textbook reference for both.
For example, PRM works as follows. A random uniform sample $\Omega$ is drawn from $E^n$ by randomly choosing values of the $n$ coordinates.  To have the sample restricted to $M$, one runs a check if the generated point is feasible, in the sense that it belongs to $M$.  So $D = \Omega \cap M$ serves as a discrete approximation of $M$.  A graph approximation $V^{(1)}_{\varepsilon}$ to $M$ may be formed from $D$ by connecting pairs of vertices according to a so-called \textit{local planner}.  For example, whenever the distance between them is bounded by a specified threshold $\varepsilon \ge 0$ and after checking that the entire geodesic segment is contained in $M$.  All of this is usually implemented to be performed at the same time as new samples in $\Omega$ are drawn.  

To solve more subtle problems than connectivity and to set up the use of discrete Morse theory, we need the following higher-dimensional generalization of the graph model.

\begin{DefRefName}{VR}{Vietoris-Rips complex}
Given a finite metric space $D$ with metric $\ell$ and a numerical threshold $\varepsilon \ge 0$, the \textit{(modified) Vietoris-Rips complex} denoted $V_{\varepsilon} (D)$ is a simplicial complex with an $n$-simplex corresponding to each $(n+1)$-tuple $S = \{ x_0, \dots, x_n \}$ whenever $\diam (S) < \varepsilon$.  For the purposes of using this construction in the same manner as PRM, we will assume that the edges of $S$ have to also pass the feasibility check.
So the 1-dimensional skeleton of $V_{\varepsilon} (D)$ is the graph approximation mentioned above, and $V_{\varepsilon} (D)$ is the clique complex completion of its 1-dimensional skeleton.
\end{DefRefName}

\begin{RemRef}{methods}
Other than building this map, PRMs are multi-query planners in that they create a library of feasible paths concurrently with sampling.  This library can be reused for solving new problems within the same environment.  RRT is a single-query planner which builds a feasible tree with the root in the initial configuration.  It stops when the last sampled configuration is within a feasible edge from the terminal point.  There are various enhanced algorithms such as RRT$^\ast$, PRM$^\ast$, Hybrid RRT-PRM, and others.  We mention them because these are the methods we will compete against in applications.  Also, we can use the implementation of the map building procedures from PRM to construct the (modified) Vietoris-Rips complex $V_{\varepsilon} (D)$.
\end{RemRef}
 
The area of applied topology called Topological Data Analysis (TDA) uses constructions similar to Vietoris-Rips complexes and the evolution of their properties across the values of the parameter $\varepsilon$ to create ``topological summaries'' reflecting the geometry of $D$ and, in an approximate way, of $M$.  Parts of the literature employing TDA flavored tools such as persistent homology have ($\dagger$) topologists addressing robotics problems as well as ($\ddagger$) roboticists using TDA.  Without making this a focus of the survey, we will mention the most prominent trends in this area.

Rob Ghrist and his students and collaborators \cite{sB:19,sBrG:15,sBrGvK:15,rG:09,rGvP:12,WSB:23} are examples of the former group ($\dagger$). There is a lot of variety to this work, but the bulk is addressing modeling motion paths for multiple robots in graphs.  This may seem a special case.  However, many standard models in robotics are based on grids and graphs $V_{\varepsilon}$ arising from sampling inside the C-space.  Also, given that many applications require safety guarantees in path planning, we will see below that one is often led to restrict planning to so-called medial axes of the C-space.  Once restricted to the graph formed by medial axes, the planning can be done entirely in the medial graph.

Florian Pokorny, Danica Kragic and others \cite{BPKE:015,CVPK:19,MPNG,OPT:21,PHR:14,PHR:15,PK:15} from their labs in KTH Royal Institute of Technology, Stockholm, are examples of the latter group ($\ddagger$).  For example, \cite{MPNG} used persistent homology in efficient computations of rigid configurations of obstacles that form energy-bounded cages of an object for grip design and other applications.  Almost universally these applications are to the sampling-based algorithms.  It is interesting to observe, and matches experience in other applications of persistent homology, that while the agorithm in loc. cit. runs in $O(s^3 + sn^2)$ time for $s$ the number of samples and $n$ the total number of vertices used to describe  the object and obstacles (so typically $n \ll s$), the observed runtimes are closer to $O(s)$ for a fixed $n$ which makes them very competitive against benchmarks.  Similarly \cite{PHR:15}, building on \cite{PK:15}, uses homology to generalize checking for connectedness properties using graph-based methods to more sophisticated classification of trajectories which includes the higher homotopy information.

One of the main current obstacles to rigorous algorithms and theorems starting from a uniform sample $D$ from $M$ is poor understanding of how to translate topological properties of $V_{\varepsilon} (D)$ to those of $M$.  Some statements are very believable and even statistically almost always true and useful in specific applications, without rigorous theory behind them.  As one striking example, consider the open question whether $V_{\varepsilon} (G)$ is contractible for large enough $\varepsilon$, where $G$ is the set of vertices in a regular cubical grid in $E^n$.  This is in contrast to the set of vertices of a Gromov hyperbolic complex which is known to be contractible for large enough $\varepsilon$.

\SecRef{Discrete Morse theory methods}{DMTM}

\SSecRef{Some preliminaries}{BDMT}
Discrete Morse theory was invented by Forman as a combinatorial version of the smooth Morse theory \cite{rF:98}.  The parallel development due to Bestvina and Brady is now understood to be essentially equivalent \cite{mZ:22} but Forman's framework is much more in use in TDA literature.  One way to generate a so-called discrete gradient vector field on a simplicial or regular cellular complex is by writing down a discrete Morse function.  It's also possible to introduce and work directly with a gradient vector field.  The geometric rewards are very efficient algorithms for computing critical cells and cell decompositions with much fewer cells than the original triangulation of the C-space. In addition to the efficiency, there is a guarantee of invariance of the homotopy type of the complexes under these reductions, so one has cellular homology computations that are in fact homology computations for the C-space. 

We will want to create simplicial or cellular analogues of the methods from motion planning that have worked well in the past.  
To start, one needs to introduce a triangulation or a cellular structure on the C-space.  Let's assume the C-space is well-sampled, so we start in our application with a Vietoris-Rips complex $\K$ of a sample from the C-space viewed as a metric subspace of a Euclidean space of certain appropriate dimension.  

The following definitions are due to Forman \cite{rF:98,Forman_guide}.

\begin{DefRef}{NMU}
A \textit{discrete Morse function} on $K$ is a real-valued 
function $f$ on the simplices of $K$
satisfying for all simplices
$\sigma$ of dimension $p$,
\begin{itemize}
\item $\displaystyle \# \{\beta^{(p+1)} > \sigma | f(\beta) \leq f(\sigma)\} \leq 1;$
\item 
$\displaystyle  \# \{\alpha^{(p-1)} < \sigma | f(\sigma) \leq f(\alpha)\} \leq 1,$
\end{itemize}
where $\#$ means the cardinality of a set, and
the notations $\alpha < \sigma$  and  $\beta  >  \sigma$ mean  
that $\alpha$ is a face of $\sigma$ and $\sigma$ is a face of $\beta$,
respectively.
We use the notation $\sigma^{(p)}$ to indicate that the dimension 
of simplex $\sigma$ is $p$.
\end{DefRef}

\begin{DefRef}{DHR}
For a discrete Morse function $f$
on $K$,
a simplex
$\sigma^{(p)}$ is called a
\textit{critical simplex}
if
\begin{itemize}
\item 
$\displaystyle\# \{\beta^{(p+1)} > \sigma \, | \, f(\beta) \leq f(\sigma)\} =0;$
\item 
$\displaystyle\# \{\alpha^{(p-1)} < \sigma \, | \, f(\sigma) \leq f(\alpha)\} =0.$
\end{itemize}
If $\sigma$ is a critical simplex,
then we call $f(\sigma)$ a \textit{critical value} of $f$.
Conversely, if a simplex $\sigma$ is not critical, 
then
we say that $\sigma$ is a \textit{regular simplex},
and $f(\sigma)$ is a \textit{regular value}.
\end{DefRef}

It turns out regular simplices occur in pairs $\alpha < \beta$ where $\alpha$ is a codimention-1 face of $\beta$ with the values of $f$ on the pair satisfying the inequality $f(\alpha) \le f(\beta)$.

Forman defines an analogue of a vector field on a complex.

\begin{DefRef}{DVF}
		A discrete vector field $\V$ on $\K$ is a collection of ordered pairs of simplices of the form $(\alpha,\beta)$ such that 
		\begin{itemize}
			\item $\alpha < \beta$,
			\item $\dim (\beta) = \dim (\alpha) +1$,
			\item each simplex from $\K$ is in at most one pair of $\V$.
		\end{itemize}
  
The flow-lines associated with a smooth vector field have a counterpart in the discrete setting we will refer to as a $\V$-path.
Given a discrete vector field $\V$, a \textit{$\V$-path} is a sequence of cells:
$$\alpha_0^{(p)}, \beta_0^{(p+1)},  \alpha_1^{(p)}, \beta_1^{(p+1)},  \alpha_2^{(p)}, \ldots, \beta_{r-1}^{(p+1)}, \alpha_r^{(p)},    
$$
where $(\alpha_i, \beta_i) \in V$, $\beta_{i}>\alpha_{i+1}$, and $\alpha_{i} \neq \alpha_{i+1}$ for all $i = 0,\ldots, r-1$.
A $\V$-path is a \emph{non-trivial closed $\V$-path} if $\alpha_r = \alpha_0$ for $r > 1$. 
A \textit{gradient vector field} on a simplicial complex $\K$ is a discrete vector field $\V$ on $\K$ which does not admit non-trivial closed $\V$-paths.
\end{DefRef}

From above we know that to a discrete Morse function $f$ there is an associated gradient vector field.

\begin{ThmRefName}{For39}{Forman, Theorem 3.9 \cite{rF:98}}
A discrete vector field is the gradient vector field of a discrete Morse function if and only if there are no non-trivial closed $\V$-paths.
\end{ThmRefName}

\begin{ThmRefName}{For34}{Forman, Theorem 3.4 \cite{rF:98}}
Suppose $\V$ is the gradient vector field of a discrete Morse function $f$. Then a sequence of simplices is a $\V$-path if and only if 
$\alpha_i < \beta_i > \alpha_{i+1}$ for $i=0, \ldots, r$, and
\[
f(\alpha_0) \ge f(\beta_0) > f(\alpha_1) \ge f(\beta_1) > \cdots \ge f(\beta_r) > f(\alpha_{r+1}).
\]
\end{ThmRefName}

\begin{DefRef}{gradpath}
If $\V$ is a gradient vector field, a $\V$-path is called a \emph{gradient path} or a \emph{trajectory}.
\end{DefRef}

It will be convenient to exploit the essential equivalence between presenting a discrete Morse function on a cellular complex $K$ and presenting its gradient vector field.  Another equivalent data containing the information about the gradient vector field is called a modified Hasse diagram.  The (unmodified) Hasse diagram is the directed graph whose vertices are the simplices of $K$ and edges are the directed codimension-one face relations in $K$, with the canonical direction from the higher dimensional simplex $\beta$ to its lower dimensional codimension-one face $\alpha$.  If $K$ comes with a discrete vector field $\V$, the modified Hasse diagram has an edge between $\alpha$ and $\beta$ assigned the opposite direction, from $\alpha$ to $\beta$, compared to the original Hasse diagram whenever the pair $(\alpha,\beta)$ is in $\V$.

\begin{ThmRefName}{For62}{Forman, Theorem 6.2 \cite{rF:98}}
		There are non-trivial closed $\V$-paths if and only if there are no
nontrivial closed directed paths in the modified Hasse diagram.
\end{ThmRefName}

\SSecRef{Relations between the smooth and discrete Morse theories}{Rel}
Of course, a theory that is discrete from the outset is much more amenable to computations in an applied field such as robotics.  The alternative is imposing smooth conditions on the environment, requiring explicit formulas for obstacle boundaries, and, further, approximating all functions by special functions such as polynomials for effective computation. 

The classical Morse theory that inspired Forman can be thought of as a way to use generic smooth functions and gradient fields on manifolds to deduce properties of the manifold.  He proved a number of fundamental theorems of discrete Morse theory in \cite{rF:98}.

Let $\V$ be a gradient vector field on a simplicial complex $\K$.  Let's include the empty set as a $(-1)$-dimensional face of any simplex, as is conventional in combinatorial literature, including some quoted references in this paper.  From inspecting the definition, a nonempty simplex $\alpha$ is a critical simplex with respect to $\V$ if one of the following holds:
	\begin{itemize}
		\item $\alpha$ does not appear in any pair of $\V$, or
		\item $\alpha$ is a $0$-simplex and $(\emptyset,\alpha) \in \V$.
	\end{itemize}

Recall that a cell is a space homeomorphic to a closed Euclidean ball.  A \emph{cellular complex} is a topological space with the quotient topology built by gluing on cells of increasing dimension to subject to some local finiteness conditions. 

\begin{ThmRefName}{For}{Forman, Theorem 2.5 \cite{rF:98}}
		If $\K$ is a simplicial complex and $\V$ is a gradient vector field on $\K$, then $\K$ is homotopy equivalent via a sequence of elementary cellular collapses to a CW complex with exactly one cell of dimension $p$ for each critical (with respect to $\V$) simplex of dimension $p$.
\end{ThmRefName}

There is a number of implementations of this algorithm together with subsequent homology computations which just as in classical Morse theory are simpler because the dimensions of cellular chain modules are greatly reduced.  In fact, this pre-processing step implemented using \texttt{MorseReduce} from Mischaikow and Nanda \cite{MiNa} is what makes persistent homology computations using \texttt{Perseus} \cite{Pers} at all feasible for large data sets.  

Several formal transitions from smooth Morse functions and fields to discrete analogues can be found in the literature.  Approaches of Benedetti \cite{bB:12,bB:16} and Mrozek et al. \cite{BKMW,MKW:16,MW:20} define them on a given triangulation of a smooth manifold.  Gallais \cite{eG:10} creates a triangulation along the flow lines of the field and defines a function in a somewhat dynamic fashion.  A bonus in the latter result is that stable manifolds are subcomplexes of triangulated neighborhoods of critical points.  This allows to translate all applied smooth methods and results from Section \refSS{PVF} to the discrete setting.  

\SSecRef{Advantages of data-driven applications of discrete Morse theory}{AdvDMT}
In addition to computability from get-go, there are advantages to starting with just a discrete metric space $D$ sampled from $X$:
\begin{enumerate}
    \item no need to assume the C-space $X$ is a manifold,
    \item the ease and flexibility of generating discrete Morse functions on complexes from the given function on the vertices $D$,
    \item the extensive library of geometrically meaningful functions on $D$ which have no computable smooth analogues.
\end{enumerate}
  
Let us address the items in order.

\smallskip

(1) It is known from Kapovich and Millson \cite{mKjM:02} that each closed smooth manifold can be generated as a component of the configuration space of a planar linkage.  On the other hand, we don't expect a generic configuration space be modeled by a single manifold.  A much more realistic model is a stratified space with manifold strata of various dimensions.  Development of methods that address global planning problems with no assumptions on local structure of $X$ gives a massive advantage.

\smallskip

(2) The discrete Morse functions can be generated on the triangulation of a manifold $M$ with $D$ as vertices or on the Vietoris-Rips complex $V_{\varepsilon} (D)$ where the vertices $D$ is a random sample from $X$.  One general algorithm \texttt{Extract} is due to King, Knudson, and Mramor.

\begin{ThmRefName}{KKM}{Infinite persistence version, \cite{hKkKmM:05}}
There is an algorithm \texttt{Extract} that
takes as input a finite simplicial complex $K$ and an injective function 
$h \colon K^{(0)} \to [0,1]$ on the vertices of $K$. Its output is three lists $A$, $B$, $C$ of simplices of $K$ so that each simplex of $K$ is on exactly one list, together with a one-to-one onto map $r \colon B \to A$ so that $r(\alpha)$ is a codimension-one face of $\alpha$.  This output encodes a modified Hasse diagram with no directed loops, so by Theorem \refT{For62} the diagram is equivalent to defining a gradient vector field on $K$, and further equivalent to defining a Morse function by Theorem \refT{For39}.  
\end{ThmRefName}

In the gradient vector field, the simplices from $C$ are critical, and the function $r$ gives a pairing of the noncritical simplices which can be used for performing a sequence of elementary collapses across each simplex $\alpha$ in $B$.  The discrete Morse function $h'$ has a larger value on
$r(\alpha)$ than on $\alpha$ and smaller on all other codimension-one faces of $\alpha$.

\begin{RemRef}{Remark}
    It will be important for us to point out that \texttt{Extract} has two stages, \texttt{ExtractRaw} and \texttt{Cancel}.  According to \cite[Theorem 3.1]{hKkKmM:05}, everything said about \texttt{Extract} in Theorem \refT{KKM} applies to \texttt{ExtractRaw} as well.  The second stage \texttt{Cancel} is a simplification stage where some critical cells are paired and cancelled in pairs.  Some facts are proved about the outcome of entire \texttt{Extract}.  For example \cite[Lemma 3.5]{hKkKmM:05} guarantees that the vertex with the minimum value of $h$ will be critical for the Morse extension $h' \colon K \to [0,1]$.  However, sometimes it's important to preserve the geometrically important information contained in the intermediate construction of $h'$ after performing just \texttt{ExtractRaw}.  One useful feature is that \texttt{ExtractRaw} makes each vertex with a local minimum critical.

    There are simplified versions of the original algorithm with significant improvements in runtime, specifically of the stage \texttt{ExtractRaw}, see Fasy et al. \cite{Fasy}.
\end{RemRef}

Another algorithm for extending a function given on vertices to a discrete Morse function on the entire complex is \texttt{ProcessLowerStars} from Wood, Robins, and Sheppard \cite{WRS:11}.  
It is valid for functions on 2- and 3-dimensional cubical complexes. We will discuss \texttt{ProcessLowerStars} in Section \refSS{Robins}.

\smallskip

(3) There are many useful functions that are intrinsic to discrete metric spaces $D$.  Of course, lots of commonly used functions on $D$ are simply restrictions of functions that can be defined on continuous models from where $D$ is sampled.  More interesting functions are the ones that require counting or computing a finite number of values from $D$.  

One example comes from \textit{kernel density estimation}.  We can start with a very natural function
\[
\delta_{\varepsilon} (x) = \# \{ x' \in D \ \vert \  d(x',x) \le \varepsilon \} 
\]
This density function by itself can serve as the function $h$.  There are more useful modifications built out of $\delta_{\varepsilon}$ for various purposes, as explained in Section \refSS{Upad}.  A more sophisticated point of view is to view them as density estimators and feed into a kernel density estimator which is dependent on parameters with statistical significance. We have had success with using these functions for creating goal posts in $X$ as critical points near boundaries of obstacles to plan efficient motion in Upadhyay et al. \cite{UGE:21,UGWE:20}.
The key idea that illustrates the advantage of this data-driven application but is likely true for other applications can be stated as follows.  While in using sampling from $X$ we no longer have formulas describing the parameter space of the environment, for example the obstacle boundaries, the density estimator allows to detect proximity to the boundary with numerous advantages such as an intrinsically local computation, approximate distance from the closest obstacle with guarantees, and others.  

A different function that is well-suited for maximizing clearance from obstacles or for designing local minimum trap avoidance is \textit{eccentricity}.  In its simplest form, an eccentricity function within a sample $D$ is given as \[ e_p(x)= \frac{1}{\# D} \sum_{x' \in D} d(x,x')^p, \quad x' \in D.\]  

Other choices for $h$ may come from the standard set of lens functions used in the Mapper algorithm \cite{CMS} in TDA. see also \cite{CMO,NLC:11,N:15}.  For example, there are scaled projections based on principal component analysis normally used for dimension reduction, or the centrality function with the effect similar to that of eccentricity.

\SecRef{Applications of discrete Morse theory}{ADMT}

Behind this section is the belief that the fundamental questions in motion planning can be addressed in a systematic ways using discrete Morse theory.  It has recently seen an explosion of fantastic applications in machine learning and topological data analysis and beginning to have the same kind of impact in robotics.  We have already argued in Section \refSS{Rel} that it provides a faithful discrete analogue of existing popular smooth methods based on artificial potential functions, with faster run speed and using less memory.  Here we survey recent genuinely novel applications.

\SSecRef{Skeletonization, sphere percolation (Robins et al.)}{Robins}
We start with the work of Vanessa Robins and collaborators \cite{RDS:15,WRS:11} that is not expressly to motion planning but to related issues in computer vision and image analysis.  It has had great influence on the use of discrete Morse theory and, as a paradigm, its applications to robotics. 

The algorithm \texttt{ProcessLowerStars} introduced in \cite{WRS:11} produces a discrete vector field $\V (K)$ starting with an arbitrary positive bounded function (a.k.a grayscale function) on cells (pixels or voxels) of a 2D or 3D pixel grid $K$.  The goal is to be able to reduce the pixel-wise information contained in the grayscale image into the meaningful layering of larger regions with grayscale value bounded from below.  These is the most meaningful visual information contained in the image.  The field $\V (K)$ is leveraged for faster and computationally less expensive extraction of TDA invariants such as persistent homology.  

The grayscale images are modeled as cubical 2D and 3D complexes with real-valued functions defined on their vertices. The pixels or voxels are made correspond to 0-cells of the cubical complex $K$.  This gives us a positive bounded function $g$ on vertices of a cubulated plane or space.  To apply the algorithm, it is necessary to have the values of $g$ on all of the 0-cells be distinct.  If they are not, one can use a linear ramp applied to $g$ to ensure they are unique.  The algorithm inspects all of the cells in the lower star of a 0-cell $x$. The lower star $L(x)$ contains all cells $\sigma \in K $ in the cubical complex such that $g(x)=\max_{y \in \sigma} g(y)$. To give an ordering to higher dimensional cells, a new function $G$ is introduced.  If $\sigma$ contains the vertices $\{x,y_{1}, \cdots,y_{n}\}$, then 
 \[
 G(\sigma)=\{g(x),g(y_{i_{1}}), \cdots, g(y_{i_{n}})\}, 
 \]
 where $g(x)>g(y_{i_{1}})>\cdots>g(y_{i_{n}})$. This will allow us to impose the lexicographical ordering on values of $G$ when performing this algorithm.
 
 The algorithm will then take all of the cells of in $L(x)$ and either pair a cell $\tau$ with a codimension-one face $\sigma$ and place them as a pair $(\sigma, \tau)$ in $\V$ or take just a single cell $\sigma$ and place it in $C$.  This is decided in the following way.
 
\begin{itemize}
    \item If $L(x)=\{x\}$, then $x \in C$.  Otherwise, take the minimal (with respect to $G$) $1$-cell, $\tau$, and $(x, \tau) \in V$. All other $1$-cells are added to a queue called PQzero, since they have no remaining unpaired faces in $L(x)$.  All cofaces of codmension 1 of $\tau$ in $L(x)$ that have 1 unpaired face in $L(x)$ are added to a different queue called PQone.
    \item Next the algorithm takes the minimal cell (with respect to $G$) that is in PQone and  either moves it to PQzero, if it has no unpaired faces remaining in $L(x)$, or if it still has an unpaired face it gets paired with that face and is put into $V$. Then it considers all cofaces of codimension-one of both cells that were just paired and put into $V$.  If any of these cofaces have exactly one unpaired face, they are added to the queue PQone.
    \item  Then, if PQzero is not empty, the algorithm takes the minimal cell (with respect to $G$), let's say $\delta$, and places it singly in $C$. Then all cofaces of codimension-one of $\delta$ are inspected.  If any of these cofaces has exactly one unpaired face, it is placed in PQone. 
\end{itemize}
The sorting loop keeps running until both PQzero and PQone are empty.  

The outcome is a discrete vector field $\V$ on the cubical complex.  Notice that the decisions in the algorithm are based entirely on values of $G$ on cells in the links of the pair in question.  This is a very local construction, so the algorithm can be applied to the interior of any subcomplex of a cubulated rectangle or rectangular prism.  This way the algorithm is applicable to any grid-based sample from a 2- or 3-dimensional C-space.

\medskip

The main results of \cite{WRS:11} verify that the collapses associated with $\V$ and the reduced cellular structure from Theorem \refT{For} respect the filtration given by the grayscale function, and so the persistent homology computation for the levels sets of that filtration can be performed on the simplified cellular level.

A later paper Delgado-Friedrichs, Robins, and Sheppard \cite{RDS:15} uses a variant of \texttt{ProcessLowerStars} to extract skeletons and partitions of images.  The applications are to analyze micro-CT images of porous materials such as spherical bead packings, sandstones, and limestone. What the authors call a skeleton is a geometric graph that passes through all gaps in the material through points that are equidistant from all closest obstacles.  Skeletonization of an image can be viewed as a 1-dimensional summary of the tunneling system embedded in the image.  The central idea is to relate persistent homology to structural measurements in porous materials such as grain and pore-size distributions, connectivity, and the critical radius of a percolating sphere. This critical radius is essential for understanding fluid flow and transport properties in porous mediums.  

The authors define in Section 3.2 \cite{RDS:15} the Morse skeleton as the subcomplex built out of unstable complexes of critical cells.  If $\alpha^{(p)}$ is a $p$-cell in $K$, they define the \textit{stable} set of $\alpha$ as the smallest set of $p$-cells $\delta^{(p)}$ in $K$ such that all $\V$-paths from $\delta$ to $\alpha$ contain all $\V$-paths leading to $\alpha$.  The the \textit{unstable} set $U_K (\alpha)$ of $\alpha$ is the dual notion.  The \textit{unstable complex} of $\alpha$ is the closure of its unstable set:
\[
W_K (\alpha) = \{\gamma \ \vert \ \gamma \le \delta \in U_K (\alpha) \}. 
\]
Now the \textit{Morse skeleton} is the union $A_K$ of all $W_K (\alpha)$ with $\alpha$ a critical cell.  According to Theorem 4 \cite{RDS:15}, there is a regular collapse from each level subcomplex in the digital image to its Morse skeleton.  This can be viewed as a ``thinning process''.  The advantage of discrete Morse theory is that this collapse does not need to be performed explicitly in order to identify the Morse skeleton.

In Section 3.3 \cite{RDS:15} the authors use the stable sets for the dual task of partitioning the vertex set of $K$ into regions which serve as ``watershed'' demarcation, by observing that each vertex is in the stable set of exactly one local minimum, which is a critical vertex $\alpha$.  This defines a partition into what they call \textit{basins}.  Basins are maximal subcomplexes of $K$ with a regular collapse onto $\alpha$.  It follows they are simply-connected and free of critical cells other than their defining $\alpha$.

\begin{figure}
(a)\includegraphics[scale=0.22]{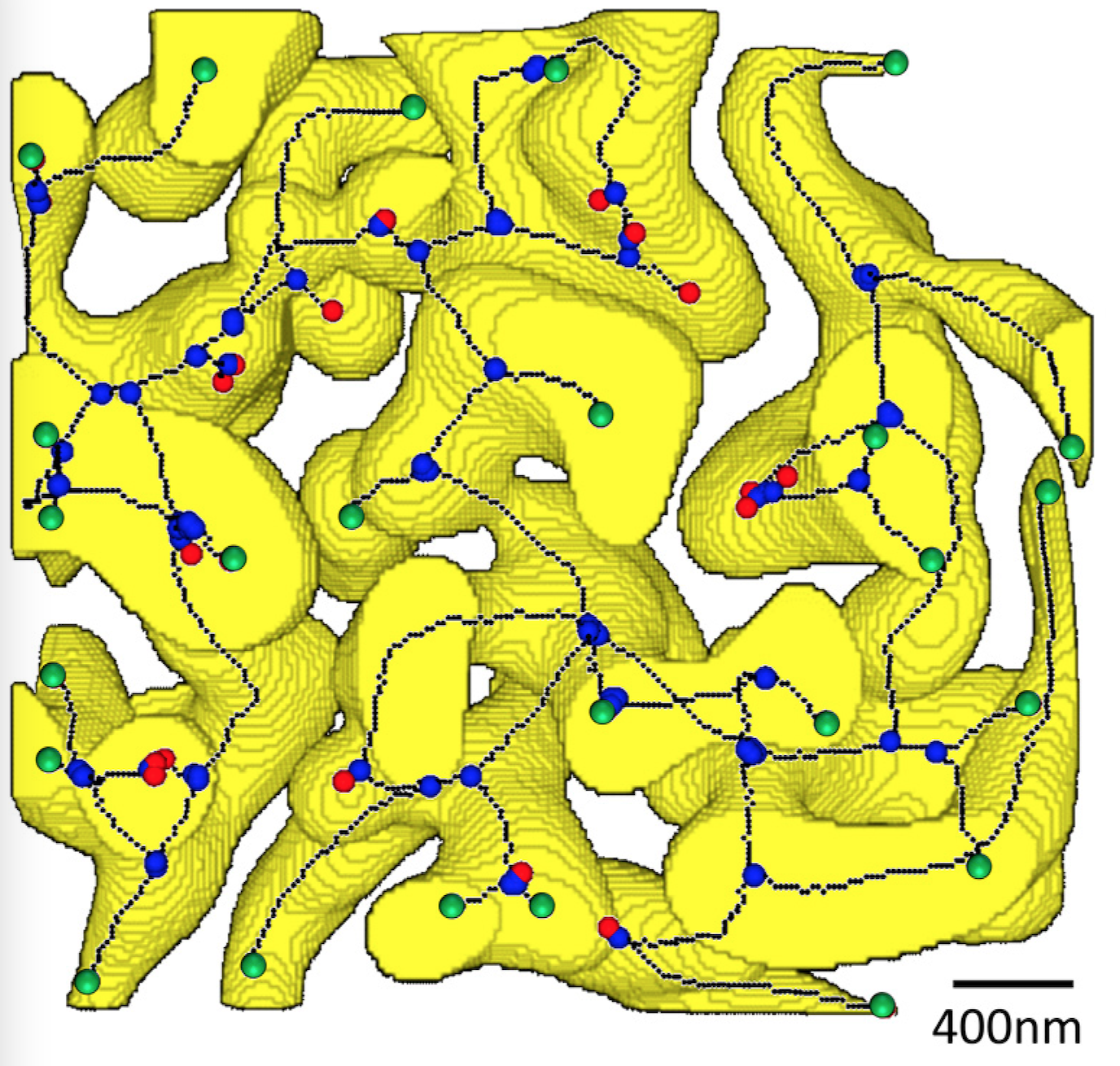}
(b)\includegraphics[scale=0.41]{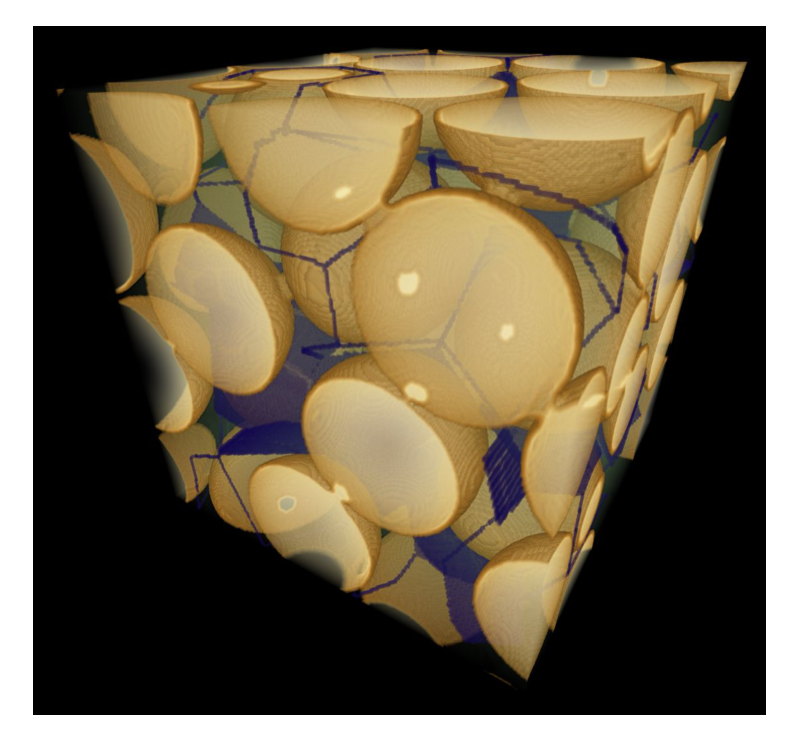}
\caption{(a) An illustration of the medial axis as a Morse skeleton and features such as junctions (blue) and ends (green and red). (b) Figure 5 from Robins et al.\,\cite{RDS:15}: a visualization of a silica 3D sphere packing with the Morse skeleton of its pore space shown in blue. }
\label{fig:envs}
\end{figure}

\begin{figure}
\includegraphics[scale=0.60]{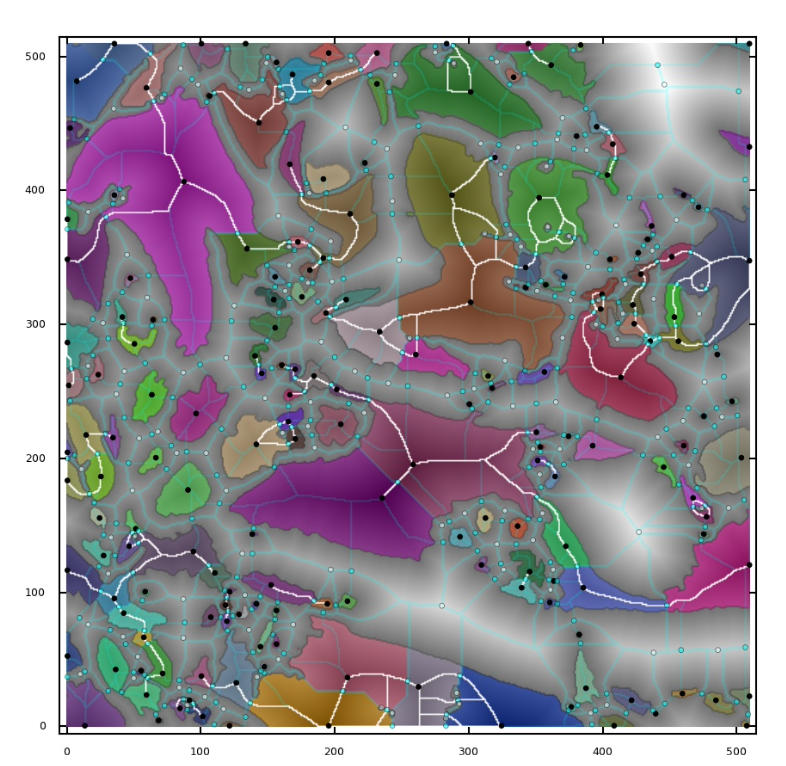}
\caption{Figure 4 from Robins et al.\,\cite{RDS:15}: illustration of the Morse analysis in a 2D Euclidean distance image (Mt. Gambier limestone); critical cells correspond to the round dots, black for minima, cyan for saddles, and white for maxima. The Morse skeleton is shown in white, and paths between critical cells in cyan. There is duality between the Morse skeleton and the cell partition determined by basins.  }
\label{fig:envs2}
\end{figure}

This setup is intimately connected to issues in robot motion path planning by interpreting the material as obstacles and the vacuous space in between as the configuration space. 
The trajectory of the center of the sphere can be interpreted as a motion path of a robot in the C-space.  In this interpretation the critical radius is the shortest possible distance from a configuration on the path to the obstacles.  This is essentially the numerical expression of the safety guarantee as a property of the motion path.  In robotics, the skeleton as defined above, is usually referred to as the \textit{medial axis}.

This work is a model solution to the medial axis problem of finding the medial axis with a subsequent solution to building a motion path restricted to the well-studied motion path problem on the medial graph.  There are two issues with this solution.  It is theoretically verified only in dimensions 2 and 3, in the sense that the algorithm spelled out above is known to generate a gradient vector field only in these dimensions.  This is discussed in the beginning of \cite{WRS:11}.  Also, it is by design useful only for the grid-type sampling method.  Of course, the other algorithm \texttt{Extract} can be applied in all dimensions and to arbitrary random samples and their modified Vietoris-Rips complexes as models of the C-space.  I don't know of any further work done in this direction.

\SSecRef{Density-based modeling (Upadhyay et al.)}{Upad}
In the sequence of papers \cite{UGE:21,UGE:22,UGWE:20,UWE:19}, C-space models are based on random samples $D$ and the resulting modified Vietoris-Rips complexes.

Let us describe a genuinely discrete planning technique which is flexible enough to achieve prescribed safety and efficiency guarantees.  It was initiated in \cite{UGE:21,UWE:19} and later perfected.  

\begin{figure}
\includegraphics[scale=1.22]{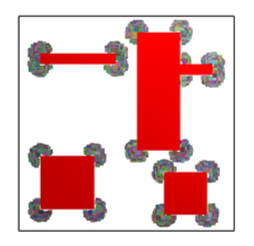}
\includegraphics[scale=1.22]{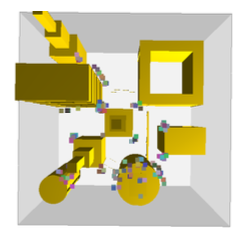}
\includegraphics[scale=0.62]{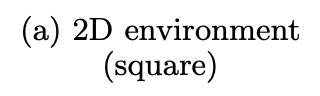}
\includegraphics[scale=0.62]{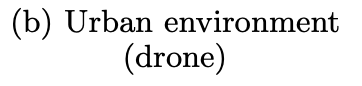}
\caption{Standard simple test environments with marked critical 0-cells to serve as goal posts for path planning.}
\label{fig:simple}
\end{figure}

The setting is a uniform sample drawn from a C-space $X$. This can be modeled on the intersection of a fine grid with the subspace $X$ of a Euclidean space $E^n$.  Assume also that the obstacles are described by linear inequalities, so that all connected components of $\mathcal{O}$ are $n$-dimensional polyhedra but don't need to be convex. Connected components of the surface of $\mathcal{O}$ are polyhedral $(n-1)$-dimensional spheres.   Observe that in this environment an optimal path will be a piece-wise linear curve with straight segments between points on the singular $(n-2)$-dimensional faces of the boundary spheres.
Observe also that the density of the sample $D$, given by the function $\delta_{\varepsilon}$ from Section \refSS{AdvDMT}, has singularities on the singular $(n-2)$-dimensional faces.  The idea is to pick out sample points in C-space that are in proximity to those vertices of polyhedral boundaries that have convex stars within $\mathcal{O}$.  Let $d(x)$ denote the shortest distance from $x$ to the closest obstacle.  Then we can use 
\[
f(x) = k {\frac{\delta_{\varepsilon}(x)}{ 1+\delta_{\varepsilon}(x)}} \exp d(x).  
\]
This function is engineered so that it has local minima at points near singular boundary points of polyhedral obstacles with convex stars.  The coefficient $k$ controls the distance to the singular points on the boundary.  These points are discovered as critical points of the Morse extension $h$ of $f$ to the Vietoris-Rips complex.  We treat these points as goal posts since they serve as breaks in the piece-wise straight trajectory.  Once we have a finite number of these points, the topological map of the environment consists of all admissible piece-wise straight paths that break at these points.  The segments used have to be admissible in the sense that they don't intersect $\mathcal{O}$.  Figure \ref{fig:simple} shows 2D and 3D environments that can be visualized, with a choice of critical 0-cells that serve as goal posts.  In the next step we let the computer sort through these paths to pick ones that satisfy imposed initial conditions.  The standard RRT- and PRM-based algorithms were applied in several simple environments.

\begin{figure}
\includegraphics[scale=0.72]{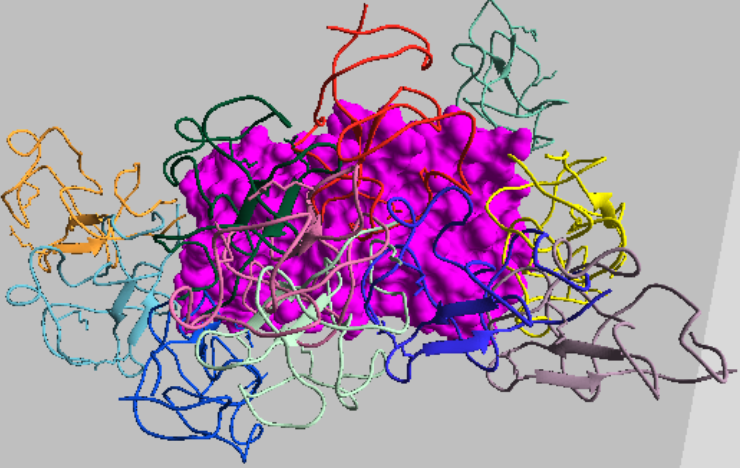}
\caption{From Figure 2 \cite{UE:22}: protein surface model is shown together with top 10 predicted 1KRN bio-molecule conformations, indicated with different colors.}
\label{fig:idp}
\end{figure}

The advantages of this simple application of the discrete Morse theory are more striking when the environments grow complex.  This happens in \cite{UE:22} where the authors generate motion paths for Intrinsically Disordered Proteins (IDPs) that need to be bound together.  Figure \ref{fig:idp} illustrates bio-molecule conformations around a protein surface model.  The problem is to find paths that implement the top 10 known conformations for a given protein.

\begin{figure}
\includegraphics[scale=0.79]{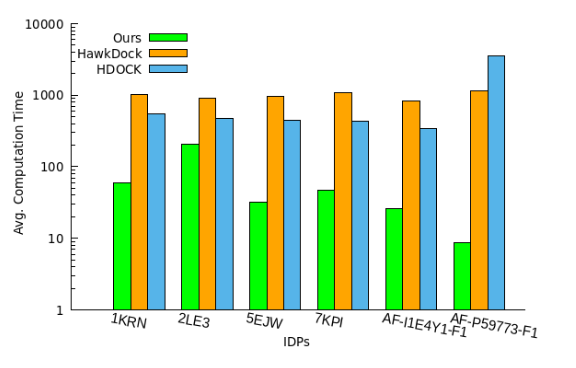} \textbf{4JUE}
\includegraphics[scale=0.79]{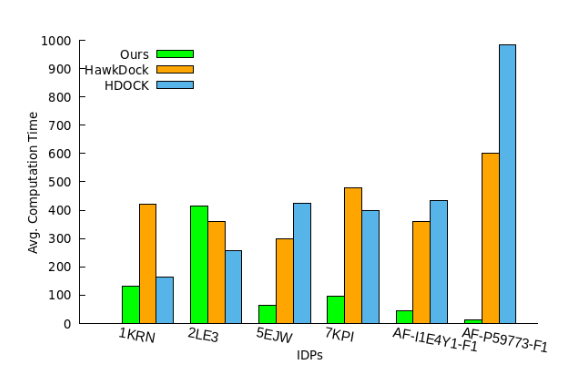} \ \ \ \ \textbf{1SQ6}
\caption{From Figure 5 \cite{UE:22}: the plots show the total computation time taken (in seconds) by three methods to predict the top 10 IDP docking conformation ensembles around the protein surface model.  The green bar shows the method with discrete Morse theory identified goal posts, the other two show results from protein-protein interaction computation servers HawkDock and HADDOCK.}
\label{fig:idpa}
\end{figure}

From the many experiments in \cite{UE:22}, we exhibit the results for two proteins, the successful but middle-of-the-road one for 4JUE and the strikingly successful one for 1SQ6 in Figure \ref{fig:idpa}.

\SSecRef{Parallelized and distributed planning}{Parallel}
C-space decompositions is an active area in robotics.  On one hand, it's the most direct way to geometrically distribute motion path planning by doing it separately in a less complex patch or cell and gluing the results together into a global solution.  The work on topological complexity is certainly related to this process.  A different point of view is what's called ``online planning'' by analogy with ``online learning" in machine learning.  It is often that the C-space is only partially known, near the initial point of the planned path.  Further information may gradually arrive through concurrent exploration or through a stream of online information.  The coverage ideas allow to start building a solution using only partial local information.

The exact cellular decomposition in Acar et al. \cite{eA:02} is based on defining a smooth function on the free space and using the levels of that function to subdivide the space into regions.  What the authors are performing in practice is a simple case of stratified Morse theory with guarantees on the topology of the regions coming from the Morse nature of the function. 
Thus simple geometry is used to compute the critical values and points, while having to deal with non-smoothness issues.  
In the beginning of Section 3.4 (Non-Smooth Boundaries) in \cite{eA:02} the authors explicitly state: ``The Morse theory that we use for decomposition is valid for the functions that are defined on a smooth compact boundary and it is not valid when the boundary is non-smooth.  However in practice, we cannot avoid non-smooth boundaries. The main problem with a non-smooth boundary is that the surface normals at the non-smooth boundary points are not defined.''  They move on to use Clarke's non-smooth gradients.   These are computationally expensive and are not well-defined in general.  We posit that these technical issues are unnecessary and are eliminated by transitioning to the discrete Morse theory as described in Section \refSS{Rel}. There are no gradients to compute, and all $\V$-paths are well-defined even in the presence of singularities.  In retrospect, the basins from Section \refSS{Robins} or other watershed decompositions such as \cite{nBgBlN:22} or \cite{FIMS} should work as well. 

\begin{figure}
\includegraphics[scale=0.92]{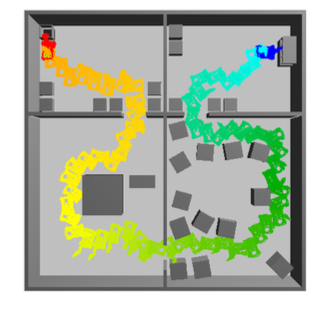}
\includegraphics[scale=0.94]{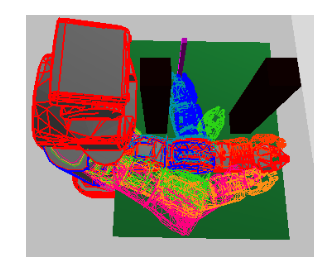}
\includegraphics[scale=0.63]{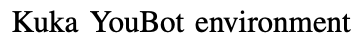} \ \ \ \ \ \ \ \ \ \ \ \ \ \ \ \ 
\includegraphics[scale=0.63]{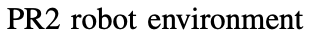}
\includegraphics[scale=0.73]{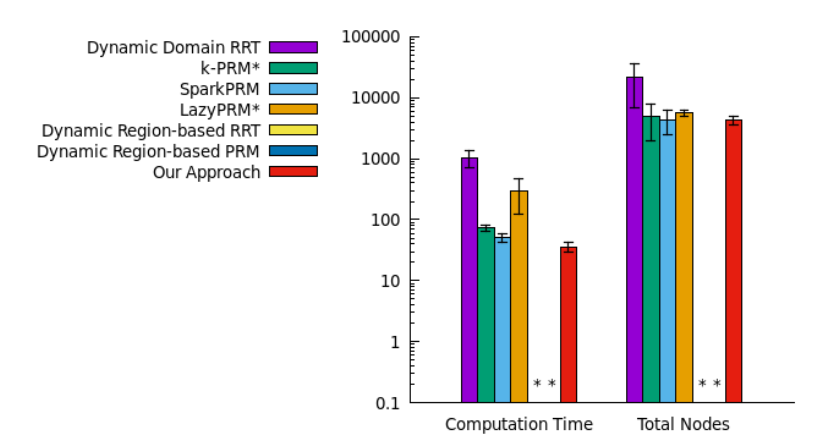}
\includegraphics[scale=0.63]{kuka-cap.png}
\includegraphics[scale=0.73]{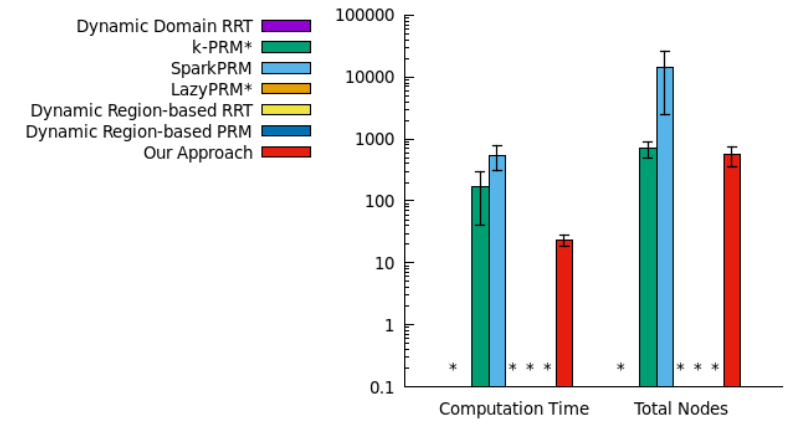}
\includegraphics[scale=0.63]{pr2-cap.png}
\caption{Incremental planning vs several state of the art algorithms.  Even in these relatively simple environments with robot a single point, the distributed planning has advantage in total computation time (which is the shortest) while the number of nodes used is the smallest.}
\label{fig:incr}
\end{figure}

Figure \ref{fig:incr} shows two examples from \cite{UGE:22}.  The coverage schemes used in \cite{UGE:22} are unfortunately still based on the generic Voronoi cells created on randomly selected nodes, but the planning was done using \cite{UGWE:20} while pasting together the Vietoris-Rips subcomplexes and the discrete vector fields generated by \texttt{ExtractRaw}.  The pasting of discrete vector fields and $\V$-paths is a procedure sorted out in \cite{dLbG:21}.

\SecRef{A collection of problems}{problems}

Some of these problems have been mentioned in the main body of the survey.  We will state the absolute versions.  This list demonstrates that the problems can range from geometric (metric) to topological (homotopy type of the configuration space) and even hybrid.  They can all be modified appropriately to ask for approximate solutions up to a given error or accuracy.
\begin{enumerate}
	\item Find the complete set of shortest paths between two configurations in C-space.  This problem can be modified by giving an allowance for paths to stray from being optimal.
	\item Find the complete set of paths up to homotopy between two configurations. In other words, compute the relative fundamental groupoid $\pi_1 (X,x_i, x_t)$ or, to lower the complexity of the problem, its abelianization $H_1 (X,x_i, x_t)$.
 Other meaningful equivalence relations on paths to consider are quantified coarse equivalences where the two paths are identified if they are fellow travelers up to a given bound.  Some definitions and observations are made in \cite{UGE:21}.
	\item Find paths that visit all marked configurations in C-space.  This problem is closely related to the Traveling Salesman Problem which has many connections with very active areas of computer science.  This circle of questions is discussed in \cite{hC:05}.
	\item \textit{The narrow passage problem.}  How does one sharpen the techniques if the most important regions  of C-space are ``narrow''?  For example, this is the problem that arises when one is designing a parallel parking algorithm for an autonomous vehicle.
	\item \textit{The safe passage problem.}  How does one modify motion path planning to include safety guarantees such as a type of buffer between the proposed path and paths that lead to collisions.  Such safety constraints can also include the restrictions of velocities when considering the underactualized movement along the path and efficiency measured in time required to complete the path.
\item Suppose the robot is a join of two subrobots, and one obtains path planning solutions for each subrobot.  How does one combine them into a motion path for the whole robot?  What properties of the path we can derive from the properties of the given solutions?  We view this as one of several ways to parallelize motion planning algorithms, which is an urgent problem in complex, multi-agent situations.
\item Suppose several robots are operating in the same environment.  How do the solutions for each robot get combined into the solution for collision-free operation of several robots?  This problem was extensively addressed in graphs by Ghrist et al. in \cite{rG:09} and more recent work.  How does one organize this problem if the robots are of different construction?  There is an interesting recent work by Orthey et al. \cite{OAT:20,OPT:21,OT:20} where the framework for multilevel motion planning is in terms of hierarchical fiber bundles.  Same kind of questions can be asked about planning around each obstacle individually and then combining the solutions into one where all obstacles are present.
\item Exact navigation systems using artificial potential functions are known to encounter the issue of getting stuck in unrecoverable state near a local singularity.  There are some well known smooth procedures to minimize occurrence of this failure.  Orthey and Toussaint \cite{OT:20} use fiber bundles mentioned in the last problem to construct what they call a local-minima tree.  This tree can be used in a failure recovery procedure by backtracking from a stalled position at one of the tree vertices.  
Is there a way to redefine the navigation function leveraging the discrete Morse theory so as to avoid stalling entirely, for example by using the eccentricity function from Section \refSS{AdvDMT}?  For this use of the eccentricity function it would be better to use the path metric in $V_1^{(1)} (D)$ to avoid sidetracking deep into the local-minima tree.

\item The main geometric questions that arise with the use of discrete samples $D$ from C-space $X$ and their Vietoris-Rips complexes are about the choice of parameters for which $X$ and $V_{\varepsilon} (D)$ have the same homotopy type.  These questions have been studied within TDA context for samples approximating manifolds embedded in a Euclidean space $E^n$, with some answers and statistical guarantees, cf. \cite{Attali,NSW}.  It's curious that there is still no complete understanding even for the regular grids in $E^n$, which are of particular importance to us in view of their use in sampling from C-space embedded in $E^n$.  Suppose $D$ is the set of integral points $\mathbb{Z}^n$ in $E^n$.  When is $V_{\varepsilon} (\mathbb{Z}^n)$ contractible?  Example 4.7 in \cite{mZ:22} provides an answer: when $\varepsilon > \sqrt{n}/(2-\sqrt{3})$.  Surprisingly, if the metric in $\mathbb{Z}^n$ is changed from Euclidean to the taxi-cab metric, which is the most natural word metric in $\mathbb{Z}^n$, the problem is open.  
\end{enumerate}

\end{document}